\documentclass[12pt]{amsart}
\usepackage{amsmath}
\usepackage{amstext}
\usepackage{amsfonts}
\usepackage{amsthm}
\usepackage{amscd}
\numberwithin{equation}{section}

\swapnumbers
\theoremstyle{plain}
\newtheorem{thm}{Theorem}[section]

\newtheorem{prop}[thm]{Proposition}
\newtheorem{cor}[thm]{Corollary}
\newtheorem{defi}[thm]{Definition}
\theoremstyle{remark}

\def\cO{{\mathcal O}}

\def\LL{\mathbb{L}}
\def\PP{\mathbb{P}}

\def\deg{\mathrm{deg} \:}

\def\dim{\mathrm{dim} \:}
\def\det{\mathrm{det} \:}

\def\im{\mathrm{im} \:}

\def\Hom{\mathrm{Hom}}
\def\Spec{\mathrm{Spec}}
\def\Ext{\mathrm{Ext}}

\def\lra{\longrightarrow}

\def\ra{\rightarrow}

\def\MM{\mathcal{M}}
\def\LL{\mathbb{L}}

\begin{document}
\title{A smooth counterexample to Nori's conjecture on the fundamental group scheme}

\author{Christian Pauly}
\address{D\'epartement de Math\'ematiques \\ Universit\'e de Montpellier II - Case Courrier 051 \\ Place Eug\`ene Bataillon \\ 34095 Montpellier Cedex 5 \\ France}
\email{pauly@math.univ-montp2.fr}

\begin{abstract}
We show that Nori's fundamental group scheme $\pi(X,x)$ does not base change correctly 
under extension of the base field for certain smooth projective ordinary curves $X$ of
genus $2$ defined over a field of characteristic $2$.
\end{abstract}

\maketitle 

\section{Introduction}

In the paper \cite{N} Madhav Nori introduced the fundamental group scheme
$\pi(X,x)$ for a reduced and connected scheme $X$ defined over an algebraically
closed field $k$ as the 
Tannaka dual group of the Tannakian category of essentially finite vector
bundles over $X$. In characteristic zero $\pi(X,x)$ coincides with 
the \'etale fundamental group, but in positive characteristic it does not (see e.g.
\cite{MS}). By analogy with the \'etale fundamental group, Nori conjectured that
$\pi(X,x)$ base changes correctly under extension of the base field. More precisely:

\bigskip

\noindent
{\bf Nori's conjecture} (see \cite{MS} page 144 or \cite{N} page 89)
If $K$ is an algebraically closed extension of $k$, then the canonical
homomorphism
\begin{equation} \label{morh}
h_{X,K}: \pi(X_K,x) \lra \pi(X,x) \times_k \Spec(K)
\end{equation}
is an isomorphism.

\bigskip

In \cite{MS} V.B. Mehta and S. Subramanian show that Nori's
conjecture is false for a projective curve with a cuspidal singularity. In this note (Corollary 4.2) we show that certain
{\em smooth} projective ordinary curves of genus $2$ defined over a field of characteristic
$2$ also provide counterexamples to Nori's conjecture.

\bigskip

The proof has two ingredients: the first is an equivalent statement of Nori's
conjecture in terms of $F$-trivial bundles due to V.B. Mehta and S. Subramanian (see
section 2) and the second is the description of the action of the Frobenius map 
on rank-$2$ vector bundles over a smooth ordinary curve $X$ of genus $2$ defined over
a field of characteristic $2$ (see section 3). In section 4 we explicitly determine the
set of $F$-trivial bundles over $X$.

\bigskip

I would like to thank V.B. Mehta for introducing me to these questions and for
helpful discussions.

\section{Nori's conjecture and $F$-trivial bundles}

Let $X$ be a smooth projective curve defined over an algebraically closed
field $k$ of characteristic $p>0$. Let $F : X \ra X$ denote the absolute Frobenius of $X$
and $F^n$ its $n$-th iterate for some positive integer $n$.

\begin{defi}
A rank-$r$ vector bundle $E$ over $X$ is said to be $F^n$-trivial if
$$ E  \ \text{stable} \qquad \text{and} \qquad F^{n*} E \cong \cO_X^{r}.$$ 
\end{defi}

\begin{prop}[\cite{MS} Proposition 3.1]
If the canonical morphism $h_{X,K}$ \eqref{morh} is an isomorphism, then any $F^n$-trivial
vector bundle $E_K$ over $X_K:= X \times_k \Spec(K)$ is isomorphic to $E_k \otimes_k K$ 
for some $F^n$-trivial vector bundle $E_k$ over $X$.
\end{prop}

\section{The action of the Frobenius map on rank-$2$ vector bundles}

We briefly recall some results from \cite{LP1} and \cite{LP2}.

\bigskip

Let $X$ be a smooth projective ordinary curve of genus $2$ defined over an
algebraically closed field $k$ of characteristic $2$. By \cite{LP2} section 2.3
the curve $X$ equipped with a level-2 structure can be uniquely 
represented by an affine equation of the form
\begin{equation} \label{eqaff}
 y^2 + x(x+1) y = x(x+1)(ax^3  + (a+b)x^2 + cx+c),
\end{equation}
for some scalars $a,b,c \in k$.
Let $\MM_X$ denote the moduli space of $S$-equivalence classes of semistable rank-$2$ vector bundles with
trivial determinant over $X$ --- see e.g. \cite{LeP}. We identify $\MM_X$ with the projective
space $\PP^3$ (see \cite{LP1} Proposition 5.1). We denote by $V : \PP^3 \dashrightarrow
\PP^3$ the rational map induced by pull-back under the absolute Frobenius $F: X \ra
X$. There are homogeneous coordinates $(x_{00}:x_{01}:x_{10}:x_{11})$ on $\PP^3$ such 
that the equations of $V$ are given as follows (see \cite{LP2} section 5)
\begin{equation} \label{ver}
V(x_{00}:x_{01}:x_{10}:x_{11}) = ( \sqrt{abc}P^2_{00}(x): \sqrt{b} P^2_{01}(x) :
\sqrt{c} P^2_{10}(x) : \sqrt{a} P^2_{11}(x)),
\end{equation}
with
$$ P_{00}(x) = x^2_{00} + x_{01}^2 + x_{10}^2 + x_{11}^2, \qquad 
P_{10}(x) = x_{00}x_{10} + x_{01}x_{11}, $$
$$ P_{01}(x) = x_{00}x_{01} + x_{10}x_{11}, \qquad
P_{11}(x) = x_{00}x_{11} + x_{10}x_{01}. $$
Given a semistable rank-$2$ vector bundle $E$ with trivial determinant, we denote by
$[E] \in \MM_X = \PP^3$ its $S$-equivalence class. The semistable boundary of 
$\MM_X$ equals the Kummer surface $\mathrm{Kum}_X$ of $X$. 
Given a degree $0$ line bundle $N$ on $X$,
we also denote the point $[N \oplus N^{-1}] \in \PP^3$ by $N$.

\begin{prop}[\cite{LP1} Proposition 6.1 (4)]
The preimage $V^{-1}(N)$ of the point $N \in 
\mathrm{Kum}_X \subset \MM_X = \PP^3$  with coordinates $(x_{00}:x_{01}:x_{10}:x_{11})$
\begin{itemize}
\item is a projective line, if $x_{00} = 0$.

\item consists of the $4$ square-roots of $N$, if $x_{00} \not= 0$.
\end{itemize}

\end{prop}

\section{Computations}

In this section we prove the following 
\begin{prop}
Let $X = X_{a,b,c}$ be the smooth projective ordinary curve of genus $2$ given by the affine model \eqref{eqaff}.
Suppose that
\begin{equation} \label{condition}
 a^2 + b^2 + c^2 + a + c = 0.
\end{equation} 
Then there exists a nontrivial family $\mathcal{E} \ra X \times S$ parametrized
by a $1$-dimensional variety $S$ (defined over k) of $F^4$-trivial rank-$2$ vector bundles with trivial determinant over $X$. Moreover any $F^4$-trivial rank-$2$ vector bundle $E$ with trivial determinant 
appears in the family $\mathcal{E}$, i.e., is of the form 
$(\mathrm{id}_X \times s)^* \mathcal{E}$ for some $k$-valued point $s : \Spec(k) \ra S$.
\end{prop}

We therefore obtain a counterexample to Nori's conjecture.

\begin{cor}
Let $X= X_{a,b,c}$ be a curve satisfying $\eqref{condition}$. 
Then for any algebraically closed extension $K$, the morphism $h_{X,K}$ is not an isomorphism
\end{cor}

\begin{proof}
Since $S$ is $1$-dimensional, there exists a $K$-valued point $s : \Spec(K) \ra S$, which is not a 
$k$-valued point. Then the bundle 
$E_K = (\mathrm{id}_X \times s)^* \mathcal{E}$ over $X_K$ is not of the form $E_k \otimes_k K$.
Now apply Proposition 2.2.
\end{proof}

\bigskip
\noindent
{\it Proof of Proposition 4.1.}  The method of the proof is to determine explicitly 
all $F^n$-trivial rank-$2$ vector bundles $E$ over $X$ for $n=1,2,3,4$. Taking 
tensor product of $E$ with $2^{n+1}$-torsion line bundles allows us to restrict
attention to $F^n$-trivial vector bundles with trivial determinant.
 
\bigskip

We first compute the preimage under iterates of $V$ of the
point $A_0 \in \PP^3$ determined by the trivial rank-$2$ vector bundle over $X$. We 
recall (see e.g. \cite{LP1} Lemma 2.11 (i)) that the coordinates of $A_0 \in \PP^3$ in the
coordinate system $(x_{00}:x_{01}:x_{10}:x_{11})$ are $(1:0:0:0)$. It follows from
Proposition 3.1  and  equations \eqref{ver}
that $V^{-1}(A_0)$ consists of the $4$ points
\begin{equation} \label{twotor}
(1:0:0:0), \ \ (0:1:0:0), \ \ (0:0:1:0) \ \ \text{and} \ \ (0:0:0:1),
\end{equation}
which correspond to the $2$-torsion points of the Jacobian of $X$.
Abusing notation we denote by $A_1$ both the $2$-torsion line bundle
on $X$ and the point $(0:1:0:0) \in \PP^3$. 

\bigskip

Both points $A_0$ and $A_1$ correspond to $S$-equivalence classes of semistable
rank-$2$ vector bundles. The set of isomorphism classes represented by the two $S$-equivalence
classes $A_0$ and $A_1$ equal $\PP \Ext^1 (A_1,A_1) \cup \{ 0 \}$ and
$\PP \Ext^1 (\cO_X,\cO_X ) \cup \{ 0 \}$ respectively, where $0$ denotes the trivial extensions $A_1 \oplus A_1$ and $\cO_X \oplus \cO_X$. Note that the two 
cohomology spaces  $\Ext^1 (A_1,A_1)$
and $\Ext^1 (\cO_X,\cO_X )$ are canonically isomorphic to $\mathrm{H}^1(\cO_X)$. 
The pull-back by
the absolute Frobenius $F$ of $X$ induces a rational map
$$F^* :\PP \Ext^1 (A_1,A_1) \lra \PP \Ext^1 (\cO_X,\cO_X ),$$
which coincides with the projectivized $p$-linear map on the cohomology $\mathrm{H}^1(\cO_X) \ra \mathrm{H}^1(\cO_X)$
induced by  the Frobenius map $F$. Since we have assumed $X$ ordinary, this $p$-linear map is bijective. Hence we obtain that there is only one (strictly) semistable 
bundle $E$ such that $[E] = A_1$ and $F^* E \cong \cO_X^2$, namely $E = A_1 \oplus
A_1$. In particular there are no $F^1$-trivial rank-$2$
vector bundles over $X$.

\bigskip

By Proposition 3.1 and using the equations \eqref{ver}, we easily obtain that
the preimage $V^{-1}(A_1)$ is a projective line $\LL \cong \PP^1$, which passes through the
two points 
$$ (1:1:1:1) \ \ \text{and} \ \ (0:0:1:1).$$
We now determine the bundles $E$ satisfying $F^*E \cong A_1 \oplus A_1$. Given
$E$ with $[F^*E] = A_1 \in \PP^3$ we easily establish the equivalence
$$ F^*E \cong A_1 \oplus A_1 \qquad \iff \qquad \dim \Hom(F^*E, A_1) = 
\dim \Hom(E, F_*A_1) = 2.$$
Suppose that $E$ is stable and $F^* E \cong A_1 \oplus A_1$. The quadratic map
$$\det : \Hom(E, F_*A_1) \lra \Hom(\det E,\det F_*A_1) = 
H^0(\cO_X(w))$$
has nontrivial fibre over $0$, since $\dim \Hom(E, F_*A_1) = 2$. Hence
there exists a nonzero $f \in \Hom(E,F_*A_1)$ not
of maximal rank. We consider the line bundle $N = \im f \subset F_* A_1$. Since
$F_* A_1$ is stable (see \cite{LaP} Proposition 1.2), we obtain the
inequalities
$$ 0 = \mu(E) < \deg N < \frac{1}{2} = \mu(F_* A_1),$$
a contradiction. Therefore $E$ is strictly semistable and $[E] = [A_2 \oplus A_2^{-1}]$
for some $4$-torsion line bundle $A_2$ with $A_2 ^{\otimes 2} = A_1$.
The $S$-equivalence class $[A_2 \oplus A_2^{-1}]$ contains three isomorphism
classes and a standard computation shows that only the decomposable bundle
$A_2 \oplus A_2^{-1}$ is mapped by $F^*$ to $A_1 \oplus A_1$. In particular
there are no $F^2$-trivial rank-$2$ bundles.

\bigskip

We now determine the coordinates of $A_2$ by intersecting the line $\LL$, 
which can be parametrized by $(r:r:s:s)$ with $r,s \in k$, with the Kummer
surface, whose equation is (see \cite{LP2} Proposition 3.1)
$$ c(x_{00}^2 x_{10}^2 + x_{01}^2 x_{11}^2) + b(x_{00}^2 x_{01}^2 + x_{10}^2 x_{11}^2)
+ a(x_{00}^2 x_{11}^2 + x_{10}^2 x_{01}^2) + x_{00}x_{01}x_{10}x_{11} = 0.$$
The computations are straightforward and will be omitted. Let $u \in k$  be a root
of the equation
\begin{equation} \label{equ}
u^2 + u = b.
\end{equation}
Then $u+1$ is the other root. The coordinates of the two $4$-torsion line bundles
(modulo the canonical involution of the Jacobian of $X$) $A_2$ such that
$A_2^{\otimes 2} = A_1$ are
$$ (u:u: \sqrt{b}:\sqrt{b}) \qquad \text{and} \qquad  (u + 1:u + 1: \sqrt{b}:\sqrt{b}).$$
Now the equation $u=0$ (resp. $u+1 = 0$) implies by \eqref{equ} $b=0$, which is 
excluded because we have assumed $X$ smooth. So by Proposition 3.1 
the preimage $V^{-1}(A_2)$ consists of the $4$ line bundles $A_3$ such that $A_3^{\otimes 2}
= A_2$. In particular there are no $F^3$-trivial rank-$2$ bundles.

\bigskip

One easily verifies that the image under the rational map $V$ given by \eqref{ver}
of the hyperplane $x_{00} = 0$ is the quartic surface given by the equation
\begin{equation} \label{eqqua}
bx_{11}^2 x_{10}^2 + c x_{11}^2 x_{01}^2 + a x_{10}^2 x_{01}^2 + x_{00}x_{10}x_{01}x_{11}
= 0.
\end{equation}
When we replace $(x_{00}:x_{01}:x_{10}:x_{11})$ with $(u:u:\sqrt{b}:\sqrt{b})$ in 
\eqref{eqqua} we obtain the equation
\begin{equation} \label{eq1}
b^2 + u^2(1+a+c) = 0.
\end{equation} 
Similarly replacing $(x_{00}:x_{01}:x_{10}:x_{11})$ with $(u+1:u+1:\sqrt{b}:\sqrt{b})$
in \eqref{eqqua} we obtain the equation
\begin{equation} \label{eq2}
b^2 + (u^2+1)(1+a+c) = 0.
\end{equation} 
Finally the product of \eqref{eq1} with \eqref{eq2} equals  (here one uses \eqref{equ})
equation \eqref{condition} up to a factor $b^2$, which we can drop since $b \not=0$ --- note
that we have assumed $X$ smooth, hence $b \not =0$ by \cite{LP2} Lemma 2.1. To summarize we have shown that if \eqref{condition} holds, then
by Proposition 3.1 there exists an $8$-torsion line bundle $A_3$ with $A_3^{\otimes 4} = A_1$ and such that
the preimage $V^{-1}(A_3)$ is a projective line $\Delta \subset \PP^3$.

\bigskip

Consider a point $[E] \in \Delta$ away from the Kummer surface --- note that $\Delta$
is not contained in the Kummer surface $\mathrm{Kum}_X$ because its intersection is 
contained in the set of $16$-torsion points. Then $E$ is stable and 
$[F^*E] = [A_3 \oplus A_3^{-1}]$.
There are three isomorphism classes represented by the $S$-equivalence class 
$[A_3 \oplus A_3^{-1}]$, namely the trivial extension $A_3 \oplus A_3^{-1}$
and two nontrivial extensions (for the details see \cite{LP1} Remark 6.2). Since 
$E$ is invariant under the hyperelliptic involution we obtain 
$F^*E = A_3 \oplus A_3^{-1}$ and finally that $E$ is $F^4$-trivial. Hence any stable
point on $\Delta$ is $F^4$-trivial.

\bigskip

Therefore, assuming \eqref{condition}, there exists a $1$-dimensional subvariety
$\Delta_0 \subset \MM_X \setminus \mathrm{Kum}_X$ parametrizing all $F^4$-trivial
rank-$2$ bundles.  Passing to an \'etale cover $S \rightarrow \Delta_0$ ensures
existence of a ``universal" family $\mathcal{E} \ra X \times S$ and we are done.

\qed 

\bigskip

{\bf Remark.} Note that equation \eqref{condition} depends on the choice of a nontrivial
$2$-torsion line bundle $A_1$. If one chooses the $2$-torsion line bundle $(0:0:1:0)$
or $(0:0:0:1)$ --- see \eqref{twotor} --- the correponding equations are
$$ a^2+ b^2 + c^2 + a+b = 0 \qquad \text{or} \qquad a^2+ b^2 + c^2 + b+c = 0.$$

\end{document}